# DONSKER THEOREMS FOR DIFFUSIONS: NECESSARY AND SUFFICIENT CONDITIONS

By Aad van der Vaart and Harry van Zanten

*Vrije Universiteit*

We consider the empirical process $\mathbb{G}_t$ of a one-dimensional diffusion with finite speed measure, indexed by a collection of functions $\mathcal{F}$. By the central limit theorem for diffusions, the finite-dimensional distributions of $\mathbb{G}_t$ converge weakly to those of a zero-mean Gaussian random process $\mathbb{G}$. We prove that the weak convergence $\mathbb{G}_t \Rightarrow \mathbb{G}$ takes place in $\ell^\infty(\mathcal{F})$ if and only if the limit $\mathbb{G}$ exists as a tight, Borel measurable map. The proof relies on majorizing measure techniques for continuous martingales. Applications include the weak convergence of the local time density estimator and the empirical distribution function on the full state space.

**1. Introduction and main results.** Let $X$ be a diffusion process on an open interval $I = (l, r) \subseteq \mathbb{R}$, that is, a strong Markov process with continuous sample paths, taking values in $I$. Denote the corresponding laws by $\{P_x : x \in I\}$ so that $X_0 = x$ under $P_x$. Assume as usual that $X$ is regular on $I$, meaning that for all $x, y \in I$ it holds that $P_x(\tau_y < \infty) > 0$, where $\tau_y = \inf\{t : X_t = y\}$. Under this condition, the scale function $s$ and the speed measure $m$ of the diffusion $X$ are well defined. The scale function is a continuous, strictly increasing function from $I$ onto $\mathbb{R}$, which implies in particular that the diffusion is recurrent. The speed measure is a Borel measure that gives positive mass to every open interval in $I$ (cf. [9, 11, 24, 25]).

We will assume throughout that the speed measure $m$ is finite, that is, $m(I) < \infty$. We denote the normalized speed measure by $\mu = m/m(I)$, and the distribution function corresponding to $\mu$ by $F$. The finiteness of $m$ implies that the process $X$ is in fact positive recurrent, and $\mu$ is the unique invariant probability measure. Hence, by the ergodic theorem, it a.s. holds

that
$$\frac{1}{t}\int_0^t f(X_u)\,du \to \int_I f\,d\mu$$

for $f \in L^1(\mu)$. It is well known that under the stated conditions, the diffusion also obeys a central limit theorem. It states that for every function $f \in L^1(\mu)$ we have the weak convergence

(1.1) $$\sqrt{t}\bigg(\frac{1}{t}\int_0^t f(X_u)\,du - \int_I f\,d\mu\bigg) \Rightarrow N(0,\Gamma(f,f))$$

as $t \to \infty$, provided that the asymptotic variance

$$\Gamma(f,f) = 4m(I)\int_I \bigg(\int_l^x f(y)\mu(dy) - F(x)\int_I f(y)\mu(dy)\bigg)^2 ds(x)$$

is finite (see, e.g., [17]). Using the Cramér–Wold device, it is easy to obtain the multidimensional extension of this result. For every finite number of functions $f_1,\ldots,f_d \in L^1(\mu)$, we have

$$\sqrt{t}\begin{pmatrix} \frac{1}{t}\int_0^t f_1(X_u)\,du - \int_I f_1\,d\mu \\ \vdots \\ \frac{1}{t}\int_0^t f_d(X_u)\,du - \int_I f_d\,d\mu \end{pmatrix}$$
$$\Rightarrow N_d\left(\begin{pmatrix} 0 \\ \vdots \\ 0 \end{pmatrix}, \begin{pmatrix} \Gamma(f_1,f_1) & \cdots & \Gamma(f_1,f_d) \\ \vdots & \ddots & \vdots \\ \Gamma(f_d,f_1) & \cdots & \Gamma(f_d,f_d) \end{pmatrix}\right),$$

where the asymptotic covariances $\Gamma(f_i,f_j)$ are defined by

(1.2) $$\Gamma(f,g) = 4m(I)\int_I \bigg(\int_l^x f(y)\mu(dy) - F(x)\int_I f(y)\mu(dy)\bigg)$$
$$\times \bigg(\int_l^x g(y)\mu(dy) - F(x)\int_I g(y)\mu(dy)\bigg)ds(x),$$

and the variances $\Gamma(f_i,f_i)$ are assumed to be finite.

In this paper we investigate the infinite-dimensional extension of the central limit theorem for diffusions. We let the function $f$ in (1.1) vary in an infinite class of functions $\mathcal{F}$, and derive necessary and sufficient conditions under which the weak convergence takes place uniformly on $\mathcal{F}$. More precisely, let $\mathcal{F} \subseteq L^1(\mu)$ be a class of functions and define for each $t > 0$ the random map $\mathbb{G}_t$ on $\mathcal{F}$ by

(1.3) $$\mathbb{G}_t f = \sqrt{t}\bigg(\frac{1}{t}\int_0^t f(X_u)\,du - \int_I f\,d\mu\bigg).$$

The map $\mathbb{G}_t$ is called the empirical process indexed by $\mathcal{F}$. If $\sup_{f \in \mathcal{F}} \int |f| \, d\mu < \infty$, the random map $\mathbb{G}_t$ is a (not necessarily measurable) random map in the space $\ell^\infty(\mathcal{F})$ of uniformly bounded functions $z : \mathcal{F} \to \mathbb{R}$, equipped with the uniform norm $\|z\|_\infty = \sup_{f \in \mathcal{F}} |z(f)|$ [see (1.6)]. We say that the class $\mathcal{F}$ is a *Donsker class* if the random maps $\mathbb{G}_t$ converge weakly in $\ell^\infty(\mathcal{F})$ to a tight, Borel measurable random element $\mathbb{G}$ of $\ell^\infty(\mathcal{F})$.

Since weak convergence in $\ell^\infty(\mathcal{F})$ to a tight Borel measurable limit is equivalent to finite-dimensional convergence and asymptotic tightness (see [1, 6], or, e.g., [29], Theorem 1.5.4), the multidimensional central limit theorem implies that the limit $\mathbb{G}$ must be a zero-mean, Gaussian random process indexed by $\mathcal{F}$ with covariance function $\mathrm{E}\mathbb{G}f\mathbb{G}g = \Gamma(f, g)$. Hence $\mathcal{F}$ can be Donsker only if there exists a version of the Gaussian process $\mathbb{G}$ that is a tight Borel measurable map into $\ell^\infty(\mathcal{F})$. By general results on Gaussian processes this is equivalent to existence of a version of $\mathbb{G}$ whose sample paths are uniformly bounded and uniformly continuous on $\mathcal{F}$ relative to the natural pseudo-metric $d_\mathbb{G}$ that $\mathbb{G}$ induces on $\mathcal{F}$, given by

$$d_\mathbb{G}^2(f, g) = \mathrm{E}(\mathbb{G}f - \mathbb{G}g)^2.$$

(In other words, the class $\mathcal{F}$ is a GC-set in the appropriate Hilbert space, in the sense of Dudley [5] or [6]. Also cf. [8], or Example 1.5.10 of [29].) Surprisingly, the existence of the limit process is also sufficient for $\mathcal{F}$ to be Donsker. In contrast with the situation for i.i.d. random elements no additional (entropy) conditions that limit the size of the class $\mathcal{F}$ are required.

It also turns out that the processes $\mathbb{G}_t$ themselves possess bounded and $d_\mathbb{G}$-continuous sample paths as well, whence the weak convergence actually takes place in the space $C_b(\mathcal{F}, d_\mathbb{G})$ of bounded, $d_\mathbb{G}$-continuous functions on $\mathcal{F}$ (cf. Theorem 1.3.10 in [29]).

THEOREM 1.1. *Suppose that $\mathcal{F}$ is bounded in $L^1(\mu)$. Then $\mathcal{F}$ is Donsker if and only if the centered, Gaussian random map $\mathbb{G}$ on $\mathcal{F}$ with covariance function $\mathrm{E}\mathbb{G}f\mathbb{G}g$ given by (1.2) admits a bounded and $d_\mathbb{G}$-uniformly continuous version. In that case, for every $x \in I$,*

$$\mathbb{G}_t \stackrel{\mathrm{P}_x}{\Longrightarrow} \mathbb{G} \qquad \textit{in } C_b(\mathcal{F}, d_\mathbb{G}) \textit{ as } t \to \infty.$$

In fact, we can prove a more general result. Since $X$ is a regular diffusion, it has continuous local time $(l_t(x) : t \geq 0, x \in I)$ with respect to the speed measure $m$. For every integrable function $f$ the occupation times formula says that

(1.4) $$\int_0^t f(X_u) \, du = \int_I f(x) l_t(x) m(dx).$$

This means that we can write the empirical process as

$$\mathbb{G}_t f = \sqrt{t} \int_I f(x)\left(\frac{1}{t}l_t(x) - \frac{1}{m(I)}\right)m(dx).$$

There is no special reason to look only at integrals of this specific type. With the same effort we can consider general integrals of the form

$$\sqrt{t}\int_I \left(\frac{1}{t}l_t(x) - \frac{1}{m(I)}\right)\lambda(dx),$$

where $\lambda$ is an arbitrary signed measure on $I$, with finite total variation $\|\lambda\|$. In this manner, we obtain a uniform central limit theorem for general additive functionals.

So let $\Lambda$ be a collection of signed measures on $I$. We define the random maps $\mathbb{H}_t$ on $\Lambda$ by

$$\mathbb{H}_t \lambda = \sqrt{t}\int_I \left(\frac{1}{t}l_t(x) - \frac{1}{m(I)}\right)\lambda(dx).$$

Slightly abusing terminology, we call $\mathbb{H}_t$ the empirical process indexed by the class $\Lambda$. By the multidimensional central limit theorem, the finite-dimensional distributions of $\mathbb{H}_t$ converge weakly to those of a Gaussian, zero-mean random map $\mathbb{H}$ on $\Lambda$ with covariance function

$$(1.5) \quad \mathrm{E}\mathbb{H}\lambda\mathbb{H}\nu = \frac{4}{m(I)}\int_I (\lambda(l,x] - F(x)\lambda(I))(\nu(l,x] - F(x)\nu(I))\,ds(x).$$

As before, the Gaussian random map $\mathbb{H}$ induces a natural pseudo-metric $d_{\mathbb{H}}^2(\lambda,\nu) = \mathrm{E}(\mathbb{H}\lambda - \mathbb{H}\nu)^2$ on the class $\Lambda$.

If the total variations of the signed measures are uniformly bounded, that is, $\sup_{\lambda \in \Lambda}\|\lambda\| < \infty$, then, for every fixed $t$,

$$(1.6) \quad \sup_{\lambda \in \Lambda}|\mathbb{H}_t\lambda| \leq \sqrt{t}\sup_{x \in I}\left|\frac{1}{t}l_t(x) - \frac{1}{m(I)}\right|\sup_{\lambda \in \Lambda}\|\lambda\| < \infty \quad \text{a.s.}$$

Hence $\mathbb{H}_t$ is a random map into the space $\ell^\infty(\Lambda)$, and we can ask whether the weak convergence of $\mathbb{H}_t$ to $\mathbb{H}$ takes place in $\ell^\infty(\Lambda)$, with a tight, Borel measurable limit process. If this is the case, we call the collection $\Lambda$ a *Donsker class*. Again, the existence of the limiting process, which is obviously necessary, is also sufficient. As before, by general results on Gaussian processes the existence can be translated into the existence of a version of the Gaussian process $\mathbb{H}$ that has bounded and $d_{\mathbb{H}}$-uniformly continuous sample paths.

THEOREM 1.2. *Suppose that $\sup_{\lambda \in \Lambda}\|\lambda\| < \infty$. Then $\Lambda$ is Donsker if and only if the centered, Gaussian random map $\mathbb{H}$ on $\Lambda$ with covariance function* (1.5) *admits a bounded and $d_{\mathbb{H}}$-uniformly continuous version. In that case, for every $x \in I$,*

$$\mathbb{H}_t \stackrel{\mathrm{P}_x}{\Longrightarrow} \mathbb{H} \qquad \text{in } C_b(\Lambda, d_{\mathbb{H}}) \text{ as } t \to \infty.$$

Theorem 1.1 is indeed a special case of Theorem 1.2, since $\mathbb{G}_t f = \mathbb{H}_t \lambda_f$, where $\lambda_f(dx) = f(x)m(dx)$.

The theory of majorizing measures provides necessary and sufficient conditions for the existence of bounded and $d_\mathbb{H}$-uniformly continuous Gaussian processes on $\Lambda$ in terms of the geometry of the pseudo-metric space $(\Lambda, d_\mathbb{H})$. See [7, 27], and Chapters 11 and 12 of [16]. We shall use this theory to prove our main theorem. Conversely, we can use it to deduce the following analytic characterization of the Donsker property.

If $(Y, d)$ is a pseudo-metric space, we denote by $B_d(y, \varepsilon)$ the closed ball around $y$ of $d$-radius $\varepsilon$.

COROLLARY 1.3. *Suppose that $\sup_{\lambda \in \Lambda} \|\lambda\| < \infty$. Then $\Lambda$ is Donsker if and only if there exists a Borel probability measure $\nu$ on $(\Lambda, d_\mathbb{H})$ such that*

$$\limsup_{\eta \downarrow 0} \sup_{\lambda \in \Lambda} \int_0^\eta \sqrt{\log \frac{1}{\nu(B_{d_\mathbb{H}}(\lambda, \varepsilon))}} \, d\varepsilon = 0.$$

PROOF. Combine Theorem 1.2 with Theorems 11.18 and 12.9 of [16]. □

In general, the majorizing measure condition is less stringent than the metric entropy condition introduced by Dudley [5]. However, the latter is often easier to work with in concrete cases. Therefore, it is useful to give a sufficient entropy condition for $\Lambda$ to be Donsker. If $(Y, d)$ is a pseudo-metric space, we denote by $N(\varepsilon, Y, d)$ the minimal number of closed balls of $d$-radius $\varepsilon$ that is needed to cover $Y$.

COROLLARY 1.4. *Suppose that $\sup_{\lambda \in \Lambda} \|\lambda\| < \infty$. Then the class $\Lambda$ is Donsker if*

$$\int_0^\infty \sqrt{\log N(\varepsilon, \Lambda, d_\mathbb{H})} \, d\varepsilon < \infty.$$

In view of definition (1.5) the covering number $N(\varepsilon, \Lambda, d_\mathbb{H})$ is the $L^2(s)$-covering number of the class of functions

$$x \mapsto \lambda(l, x] - F(x)\lambda(I), \qquad \lambda \in \Lambda.$$

These functions are of uniformly bounded variation and hence the full class, with the elements of $\Lambda$ of uniformly bounded variation, possesses a finite $L^2(Q)$-entropy integral for any finite measure $Q$. (See, e.g., [29], Theorem 2.7.5.) Unfortunately, this observation is useless in the present situation, as under our conditions the measure defined by the scale function $s$ is unbounded. Under appropriate bounds on the tails of the envelope function of the class, it is still possible to exploit the fact that the functions are of bounded variation by a partitioning argument, as in Corollary 2.7.4 of [29]. Alternatively, for special $\Lambda$ we can use the preceding corollary in combination with VC-theory. However, the best results are obtained through direct application of Theorem 1.2, as this allows to exploit the fine properties of Gaussian processes. We illustrate this in Section 2 by several examples of interest.

The Donsker theorem is based on approximation by a continuous local martingale and an analysis of local time. In Section 3 we present a uniform central limit theorem for continuous local martingales under a majorizing measure condition. This extends a result by Nishiyama [22], and is of interest on its own. In Section 4 we recall the necessary results on local time. Following these preparations the final section gives the proofs of the main results.

Let us remark that because of the use of local time, our approach is limited to the one-dimensional case. In higher dimensions one has to resort to different methods, using for instance the generator of the diffusion to relate the empirical process to a family of local martingales. This approach was followed (for general stationary, ergodic Markov processes) by Bhattacharya [3] to obtain a functional central limit theorem for

$$\frac{1}{\sqrt{n}} \int_0^{nt} \left( f(X_u) \, du - \int_I f \, d\mu \right),$$

where $f$ is one fixed function, $t \geq 0$ and $n \to \infty$. It is not clear, however, whether necessary and sufficient conditions can be obtained in this way.

The notation $a \lesssim b$ is used to denote that $a \leq Cb$ for a constant $C$ that is universal, or at least fixed in the proof.

**2. Examples.** In this section we consider four special cases of Theorem 1.2.

2.1. *Diffusion local time.* The first example is a uniform central limit theorem for diffusion local time. The space of continuous functions on a compact set $J \subseteq \mathbb{R}$, endowed with the supremum norm, is denoted by $C(J)$.

THEOREM 2.1. *Suppose that $\int_I F^2(1-F)^2 \, ds < \infty$. Then, for all $x \in I$ and compact $J \subseteq I$,*

$$\sqrt{t}\left(\frac{1}{t}l_t - \frac{1}{m(I)}\right) \xrightarrow{P_x} \mathbb{G}$$

*in $C(J)$, where $\mathbb{G}$ is a zero-mean Gaussian random map with covariance function*

$$E\mathbb{G}(x)\mathbb{G}(y) = \frac{4}{m(I)} \int_I (\mathbb{1}_{[x,r)} - F)(\mathbb{1}_{[y,r)} - F) \, ds.$$

PROOF. We apply Corollary 1.4 with $\Lambda = \{\delta_x : x \in J\}$, where $\delta_x$ is the Dirac measure concentrated at $x$. The integrability of the function $F^2(1-F)^2$ is equivalent to the finiteness of the covariance function of the limit $\mathbb{G}$. To verify the entropy condition, observe that the pseudo-metric $d$ that is induced by $\mathbb{G}$ on $\Lambda$ is given by $d(\delta_x, \delta_y) = \sqrt{|s(x) - s(y)|}$. It follows that the space $(\Lambda, d)$ is isometric to $(s(J), \sqrt{|\cdot|})$. Since $s(J)$ is compact, this implies that the entropy condition of Corollary 1.4 is satisfied. Hence, we have weak convergence in $\ell^\infty(J)$, and therefore also in $C(J)$, since diffusion local time is continuous in the space variable (see Section 4). $\square$

We remark that the weak convergence of the normalized local time process, as in the preceding theorem, cannot be extended to uniformity on the entire state space $I$. By the continuous mapping theorem, uniform weak convergence in $\ell^\infty(I)$ would imply weak convergence of the sup-norm to a finite limit. But since the function $x \mapsto l_t(x)$ vanishes outside the range of $(X_s : 0 \le s \le t)$, which is strictly within $I$ a.s., we have a.s.

$$\left\| \sqrt{t}\left(\frac{1}{t} l_t - \frac{1}{m(I)}\right) \right\|_\infty \ge \frac{\sqrt{t}}{m(I)} \to \infty,$$

which would lead to a contradiction.

On the other hand, we can construct a version of the limit process $\mathbb{G}$ with continuous (but not necessarily bounded) sample paths on the entire state space. Then the process $\sqrt{t}(l_t/t - 1/m(I))$ indexed by $I$ converges to $\mathbb{G}$ relative to the topology of uniform convergence on compacta. [To construct a version of $\mathbb{G}$ with continuous sample paths on $I$, first construct an arbitrary version $\mathbb{G}$ indexed by a countable dense subset $Q \subset I$. In view of the entropy bound obtained in the proof of Theorem 2.1 the modulus of continuity $\sup_{s,t \in J \cap Q : |s-t| < \delta} |\mathbb{G}(s) - \mathbb{G}(t)|$ of the restriction of this process to a given compact $J \subset I$ converges to zero in mean as $\delta \downarrow 0$. Thus up to a null set the process $\mathbb{G}$ is uniformly continuous on bounded subsets of its (countable, dense) index set. We can extend it by continuity to the whole state space $I$.]

For later reference we note that, given the integrability of the function $F^2(1-F)^2$, there exist positive constants $c_1, c_2$ such that, for all $x \in I$,

(2.1) $$c_1(1 + |s(x)|) \le \mathbb{E}\mathbb{G}^2(x) \le c_2(1 + |s(x)|).$$

Because the function $s$ is unbounded, this too shows that there is no version of $\mathbb{G}$ with bounded sample paths.

2.2. *Empirical process indexed by functions.* In this section we give a sufficient condition for the weak convergence of the empirical process (1.3) indexed by a general class $\mathcal{F}$ of functions. This covers many concrete examples. However, for special classes of functions, such as indicators in the line, the result can be improved, as illustrated in the next sections.

Let $(1+\sqrt{|s|})\,dF$ denote the measure with density $(1+\sqrt{|s|})$ relative to $F$.

THEOREM 2.2. *Suppose that $\int_I F^2(1-F)^2\,ds < \infty$. Then every class of functions $\mathcal{F} \subseteq L^1((1+\sqrt{s}\,)\,dF)$ that satisfies the entropy condition*

$$\int_0^\infty \sqrt{\log N(\varepsilon, \mathcal{F}, L^1((1+\sqrt{|s|}\,)\,dF))}\,d\varepsilon < \infty$$

*is Donsker.*

PROOF. In view of the occupation times formula $\mathbb{G}_t f = \sqrt{t}\int_I f(l_t/t - 1/m(I))\,dm$. Therefore, a version of the limit process $\mathbb{H}$ must be given by $\mathbb{H} f = \int f\mathbb{G}\,dm$, for $\mathbb{G}$ the limit process of the diffusion local time process obtained in Theorem 2.1. Because $\int |f|(1+\sqrt{|s|})\,dm < \infty$ by assumption and $\mathrm{E}|\mathbb{G}| \lesssim 1 + \sqrt{|s|}$ by (2.1), this process is indeed well defined. It is easily shown that this process $\mathbb{H}$ is a mean-zero process with the correct covariance structure, whence it suffices to check that it possesses a version with bounded and uniformly continuous sample paths.

Now

$$\mathrm{E}(\mathbb{H}f - \mathbb{H}g)^2 = \frac{2}{\pi}(\mathrm{E}|\mathbb{H}f - \mathbb{H}g|)^2 \lesssim \left(\int |f - g|\mathrm{E}\mathbb{G}\,dm\right)^2.$$

In view of (2.1) the intrinsic metric $d_\mathbb{H}(f,g)$ is bounded above by a multiple of the $L^1((1+\sqrt{|s|}\,)\,dm)$-norm of $f - g$. Hence the existence of the appropriate version of $\mathbb{H}$ follows from [5]. □

EXAMPLE 2.3. As a particular example, we may take any VC-class $\mathcal{F}$ with an envelope function $\mathbb{F}$ such that

$$\int_I \mathbb{F}(x)(1+\sqrt{|s(x)|}\,)\,dm(x) < \infty.$$

Then the covering number $N(\varepsilon Q\mathbb{F}, \mathcal{F}, L^1(Q))$ is bounded by $C(1/\varepsilon)^V$ for $V+1$ the VC-index of the class $\mathcal{F}$ and $C$ a constant depending on $V$ only, and any $\sigma$-finite measure $Q$ such that $Q\mathbb{F} < \infty$. (See, e.g., Theorem 2.6.7 in [29], where it is clear from the proof that the result extends to $\sigma$-finite measures $Q$.) In particular, the entropy condition of the preceding theorem is satisfied, and hence $\mathcal{F}$ is Donsker.

EXAMPLE 2.4. Another example is given by the collection of all monotone functions $f: I \to [0,1]$. Because this has a finite entropy integral for any finite measure, this class is Donsker if $\int_I \sqrt{|s|}\,dF < \infty$.

EXAMPLE 2.5. A third example is given by the collection of all functions $f : I \to [0,1]$ with $|f(x) - f(y)| \leq |x-y|^\alpha$ for some $\alpha > 1/2$ in the case that $\bar{I}$ is compact. This class has entropy relative to the uniform norm bounded above by a multiple of $(1/\varepsilon)^{1/\alpha}$ and hence satisfies the entropy condition of the preceding theorem if $\int_I \sqrt{|s|}\, dF < \infty$.

Using the approach of Corollary 2.7.4 of [29], this can be extended to unbounded state space $I = \mathbb{R}$ under the condition that for some $p < 2/3$

$$\sum_{j=1}^\infty \left( \int_{j < |x| \leq j+1} (1 + \sqrt{|s|(x)}\,)\, dF(x) \right)^p < \infty.$$

Analogy with the case of empirical processes for independent observations suggests that the class will remain Donsker if this holds for $p = 2/3$, but we have not investigated this.

2.3. *Local time density estimator.* Suppose that the invariant probability measure $\mu$ has a locally bounded density $f$ with respect to a measure $\nu$ on $I$. Then it follows from the occupation times formula (1.4) that the empirical measure $\mu_t$, defined by

$$\mu_t(B) = \frac{1}{t} \int_0^t \mathbb{1}_B(X_u)\, du,$$

has the (random) density

$$f_t(x) = \frac{m(I) f(x) l_t(x)}{t}$$

with respect to $\nu$. In the statistical literature this density $f_t$ is often called the local time estimator of $f$; see, for example, [4]. If $\nu$ is the Lebesgue measure on $I$ and $f$ is continuous, then $f_t$ is simply the derivative of the empirical distribution function.

Kutoyants [14] and Negri [19] studied the statistical properties of the local time estimator for regular diffusions on $\mathbb{R}$ that are generated by certain stochastic differential equations. In particular, for the special class of diffusions he considered, Kutoyants [14] showed that the normalized difference $\sqrt{t}(f_t - f)$ converges weakly to a Gaussian limit, uniformly on the whole state space $I$. In this section we complement and generalize their results, giving precise conditions for general regular diffusions.

The finite-dimensional distributions of $\sqrt{t}(f_t - f)$ converge weakly to those of the centered, Gaussian random map $\mathbb{H}$ with covariance function

$$\mathbb{E}\mathbb{H}(x)\mathbb{H}(y) = 4m(I) f(x) f(y) \int_I (\mathbb{1}_{[x,r)} - F)(\mathbb{1}_{[y,r)} - F)\, ds,$$

provided that these covariances are finite. The following theorem gives necessary and sufficient conditions under which this finite-dimensional convergence can be extended to uniform convergence, on compacta or on the full state space $I$. Recall that we assume throughout that $f$ is bounded on compact subsets of $I$.

THEOREM 2.6. (i) *We have $\sqrt{t}(f_t - f) \stackrel{P_x}{\Longrightarrow} \mathbb{H}$ in $\ell^\infty(J)$ for every compact $J \subseteq I$ and $x \in I$ if and only if $\int_I F^2(1-F)^2 \, ds < \infty$.*

(ii) *We have $\sqrt{t}(f_t - f) \stackrel{P_x}{\Longrightarrow} \mathbb{H}$ in $\ell^\infty(I)$ for every $x \in I$ if and only if $\mathbb{H}$ admits a version such that $\mathbb{H}(x) \to 0$ almost surely as $x \downarrow l$ or $x \uparrow r$.*

PROOF. Because $\sqrt{t}(f_t - f) = f\sqrt{t}(l_t/t - 1/m(I))$, a version of the limit process $\mathbb{H}$ of $\sqrt{t}(f_t - f)$ can be defined as $\mathbb{H} = f\mathbb{G}$, for $\mathbb{G}$ the limit process of the local time process appearing in Theorem 2.1. In the following we use a version $\mathbb{H} = f\mathbb{G}$ obtained from a version of $\mathbb{G}$ with continuous sample paths on the entire state space $I$.

For any $x \in I$,

$$E\mathbb{H}^2(x) = 4m(I)f^2(x)\left(\int_{(l,x]} F^2 \, ds + \int_{(x,r)} (1-F)^2 \, ds\right).$$

Therefore, the condition $\int_I F^2(1-F)^2 \, ds < \infty$ is equivalent to the finiteness of $E\mathbb{H}^2(x)$ for some $x$ with $f(x) > 0$ (and then for all $x \in I$), whence the condition is certainly necessary.

(i) Since $f$ is locally bounded, the map $\ell^\infty(J) \to \ell^\infty(J)$ defined by $z \mapsto fz$ is continuous for the uniform norm. Because $\mathbb{G}$ is a tight Borel measurable element in $\ell^\infty(J)$, so is the process $\mathbb{H} = f\mathbb{G}$. Thus the assertion follows from Theorem 2.1.

(ii) From (2.1) it follows that $E\mathbb{H}^2(x) = f^2(x)E\mathbb{G}^2(x)$ is bounded on $I$ only if the function $f^2s$ is bounded. Because $s(x) \to \pm\infty$ as $x$ approaches the boundary of $I$, it follows that in this case $f(x) \to 0$ at the boundary of $I$.

Because the sample paths $x \mapsto l_t(x)$ of local time vanish for $x$ near the boundary of the state space $I$ and $f(x) \to 0$ as $x$ tends to this boundary, the sample paths of the process $\sqrt{t}(f_t - f)$ tend to zero at the endpoints of $I$. If $\sqrt{t}(f_t - f)$ converges to a tight limit $\mathbb{H}$ in $\ell^\infty(I)$, then the sample paths of $\mathbb{H}$ must tend to zero at the boundary points also, as can be seen, for instance, from an almost sure construction. Thus the condition in (ii) is necessary.

To prove sufficiency, it suffices to show that there exists a version of $\mathbb{H}$ that is a tight, Borel measurable map into $\ell^\infty(I)$. Let $J_m$ be an increasing sequence of compact intervals with $J_m \uparrow I$, and let $\mathbb{H}_m = f\mathbb{G}1_{J_m}$ be the process indexed by $I$ with sample paths equal to $f\mathbb{G}$ on $J_m$ and equal to zero outside $J_m$. Because the restriction of $\mathbb{G}$ to $J_m$ is a tight, Borel measurable map into $C(J_m) \subset \ell^\infty(J_m)$ and $\mathbb{H}_m$ is the image of this restriction under the continuous map $z \mapsto fz1_{J_m}$ from $\ell^\infty(J_m)$ to $\ell^\infty(I)$, the process $\mathbb{H}_m$ is a tight, Borel measurable map into $\ell^\infty(I)$. The process $\mathbb{H} = f\mathbb{G}$ as constructed in the first part of the proof is separable, because it possesses $d_{\mathbb{H}}$-uniformly continuous sample paths on every (Euclidean) compact interval $J \subset I$, which is $d_{\mathbb{H}}$-totally bounded by tightness of $\mathbb{H}$. This implies that this version of the limit process satisfies $\sup_{x \notin J_m} |\mathbb{H}(x)| \to 0$ almost surely, as $m \to \infty$, as does the version of $\mathbb{H}$ in the statement of (ii). Consequently,

$$\sup_{x \in I} |\mathbb{H}_m(x) - \mathbb{H}(x)| \leq \sup_{x \notin J_m} |\mathbb{H}(x)| \to 0,$$

almost surely. We conclude that the process $\mathbb{H}$ is the almost sure limit in $\ell^\infty(I)$ of a sequence of tight, Borel measurable maps into $\ell^\infty(I)$. This implies that $\mathbb{H}$ is itself also a tight, Borel measurable map into $\ell^\infty(I)$, in view of the lemma below. □

The following lemma gives an easily verifiable sufficient condition for the convergence of $\sqrt{t}(f_t - f)$ on the entire state space, which is necessary under a mild regularity condition.

COROLLARY 2.7. *Suppose that $\int_I F^2(1-F)^2 \, ds < \infty$. Then if the function $f^2(x)|s(x)| \log \log |s(x)| \to 0$ as $x \to l$ or $x \to r$, the convergence*

$$\sqrt{t}(f_t - f) \stackrel{\mathrm{P}_x}{\Longrightarrow} \mathbb{H}$$

*takes place in $\ell^\infty(I)$ for every $x \in I$. If the function $f^2 s$ is monotone near $l$ and $r$, then these conditions are also necessary.*

PROOF. We first prove sufficiency. Let $\mathbb{G}$ be the limit process in Theorem 2.1, $\mathbb{H} = f\mathbb{G}$ and let $\mathbb{W}$ be a two-sided Brownian motion, emanating from zero. By the preceding theorem, it suffices to show that the sample paths of $\mathbb{H}$ converge to zero at the boundary points of $I$. Observe that $\mathrm{E}(\mathbb{G}(x) - \mathbb{G}(y))^2 = (4/m(I))\mathrm{E}(\mathbb{W}(s(x)) - \mathbb{W}(s(y)))^2$ for all $x, y \in I$. Moreover, by (2.1), we have $\mathrm{E}\mathbb{G}^2(x) \lesssim \mathrm{E}\mathbb{W}^2(s(x))$ for $x$ such that $|s(x)|$ is bounded away from 0. It follows that for $y \geq x \in I$ close enough to $r$,

$$\mathrm{E}(\mathbb{H}(x) - \mathbb{H}(y))^2$$
$$= \mathrm{E}(f(x)\mathbb{G}(x) - f(y)\mathbb{G}(y))^2$$
$$\lesssim (f(x) - f(y))^2 \mathrm{E}\mathbb{G}^2(x) + f^2(y)\mathrm{E}(\mathbb{G}(x) - \mathbb{G}(y))^2$$
$$\lesssim (f(x) - f(y)r)^2 \mathrm{E}\mathbb{W}^2(s(x)) + f^2(y)\mathrm{E}(\mathbb{W}(s(x)) - \mathbb{W}(s(y)))^2$$
$$= (f(x) - f(y))^2 \mathrm{E}\mathbb{W}^2(s(x)) + f^2(y)\mathrm{E}(\mathbb{W}(s(x)) - \mathbb{W}(s(y)))^2$$
$$\quad - 2f(y)(f(x) - f(y))\mathrm{E}\mathbb{W}(s(x))\mathbb{W}((s(y)) - \mathbb{W}(s(x)))$$
$$= \mathrm{E}(f(x)\mathbb{W}(s(x)) - f(y)\mathbb{W}(s(y)))^2,$$

by the independence of the Brownian increments. If we define $\mathbb{H}(r) = 0$ and $f(r)\mathbb{W}(s(r)) = 0$, then the processes $\mathbb{H}$ and $f\mathbb{W} \circ s$ are continuous in $L^2$ at $r$, as $s(x)f^2(x) \to 0$ as $x \uparrow r$ by assumption. Consequently, under this extension of the index set the inequality in the display remains valid for $x, y \in [x_0, r]$, for sufficiently large $x_0$. It follows that

$$\mathrm{E} \sup_{x \geq x_0} |\mathbb{H}(x)| \leq \mathrm{E} \sup_{x_0 \leq x, y \leq r} (\mathbb{H}(x) - \mathbb{H}(y))$$

$$\leq 2\mathrm{E} \sup_{x_0 \leq x \leq r} \mathbb{H}(x) \leq 2\mathrm{E} \sup_{x_0 \leq x \leq r} (f(x)\mathbb{W}(s(x))),$$

by Slepian's lemma. By the law of the iterated logarithm for Brownian motion and the condition on $f$, we have that $f(x)\mathbb{W}(s(x)) \to 0$ almost surely as $x \uparrow r$. Therefore, the median of the variables $\sup_{x_0 \leq x \leq r} |f(x)\mathbb{W}(s(x))|$ converges to zero as $x_0 \uparrow r$. As a consequence of Borell's inequality the mean of a supremum of a separable Gaussian process is bounded above by a multiple of the median and hence the right-hand side of the preceding display converges to zero. We conclude that the sequence $\sup_{x \geq x_0} |\mathbb{H}(x)|$ converges to zero in mean, and by monotonicity also almost surely, for $x \to r$. Similar reasoning applies for $x \to l$.

Now assume that $f^2(x)s(x) \downarrow 0$ as $x \uparrow r$. Then $f$ is also decreasing near $r$. Furthermore, for $y \geq x$,

$$\mathrm{E}(\mathbb{G}(y) - \mathbb{G}(x))\mathbb{G}(x) = -\frac{4}{m(I)} \int \mathbb{1}_{[x,y)}(\mathbb{1}_{[x,r)} - F)\,ds \leq 0.$$

We conclude that, for $y \geq x$ sufficiently close to $r$,

$$\mathrm{E}(\mathbb{H}(y) - \mathbb{H}(x))^2 = \mathrm{E}(f(x)\mathbb{G}(x) - f(y)\mathbb{G}(y))^2$$

$$\geq f(y)^2 \mathrm{E}(\mathbb{G}(y) - \mathbb{G}(x))^2 + (f(y) - f(x))^2 \mathrm{E}\mathbb{G}^2(x)$$

$$\gtrsim f(y)^2(s(y) - s(x)).$$

Let $x_n$ be such that $s(x_n) = e^n$. Then $s(x_n) - s(x_m) \geq s(x_n)(1 - e^{-1})$ for $n > m$, whence for sufficiently large $m$ and $n > m$,

$$\mathrm{E}(\mathbb{H}(x_n) - \mathbb{H}(x_m))^2 \gtrsim f^2(x_n)s(x_n) =: a_n^2,$$

and hence $d_{\mathbb{H}}(x_k, x_l) \gtrsim a_{2n}$ for all $n \leq k, l \leq 2n$. So the points $x_n, x_{n+1}, \ldots, x_{2n}$ are $a_{2n}$-separated, and Sudakov's inequality implies that

$$\mathrm{E} \sup_{n \leq k \leq 2n} |\mathbb{H}(x_k)| \gtrsim a_{2n}\sqrt{\log n} \gtrsim a_{2n}\sqrt{\log 2n}.$$

If $\mathbb{H}(x) \to 0$ almost surely as $x \uparrow r$, then the left-hand side tends to zero, and we conclude that $a_n^2 \log n \to 0$ as $n \to \infty$. Together with the monotonicity of $f^2s$ this implies the necessity of the right tail condition. The condition on the left tail can be seen to be necessary in the same way. □

LEMMA 2.8. *Let $X_n, X : \Omega \to \mathcal{D}$ be maps from a complete probability space $(\Omega, \mathcal{F}, \mathrm{P})$ into a complete metric space $\mathcal{D}$. If $X_n$ is Borel measurable and tight for every $n$, and $d(X_n, X) \to 0$ in outer probability, then $X$ is Borel measurable and tight.*

PROOF. The map $X$ is Borel measurable, because the convergence in outer probability implies the existence of a subsequence that converges almost surely. The pointwise limit of a sequence of Borel measurable maps into a metric space is itself Borel measurable.

If $\mathrm{P}^*(d(X_n, X) \geq \delta) \to 0$ for every $\delta > 0$, then there exists a sequence $\delta_n \downarrow 0$ such that $\mathrm{P}^*(d(X_n, X) \geq \delta_n) \to 0$. Hence given some $\varepsilon > 0$ we can find a subsequence $n_1 < n_2 < \cdots$ such that $\mathrm{P}^*(d(X_{n_j}, X) \geq \delta_{n_j}) < \varepsilon 2^{-j}$ for every $j \in \mathbb{N}$. By the tightness of $X_n$ for a fixed $n$, we can find a compact set $K_n$ with $\mathrm{P}(X_n \notin K_n) < \varepsilon 2^{-n}$.

The set $C = \bigcap_j K_{n_j}^{\delta_{n_j}}$, where $K^\delta = \{x : d(x, K) < \delta\}$, is totally bounded. If this were not the case, there would be $\eta > 0$ and a sequence $\{x_m\} \subset C$ with $d(x_m, x_{m'}) > \eta$ for every $m \neq m'$. Fix $j$ such that $4\delta_{n_j} < \eta$. There exists $\{y_m\} \subset K_{n_j}$ with $d(x_m, y_m) < \delta_{n_j}$ for every $m$, and by compactness of $K_{n_j}$ this has a converging subsequence. The tail of the sequence $x_m$ would be in a ball of radius $2\delta_{n_j}$ around the limit, which contradicts the construction of $\{x_m\}$. Thus $C$ is totally bounded, and hence its closure is compact.

If $X_n \in K_n$ for every $n$ and $d(X_{n_j}, X) < \delta_{n_j}$ for every $j$, then $X \in C$. We conclude that $\mathrm{P}(X \notin C) < 2\varepsilon$. □

2.4. *Empirical distribution function.* Let $J$ be an arbitrary subset of $I$. The empirical process $\mathbb{G}_t$ indexed by the class of functions $\mathcal{F} = \{\mathbb{1}_{(l,x]} : x \in J\}$ is the restriction of $\sqrt{t}(F_t - F)$ to $J$, where $F_t$ is the empirical distribution function, defined by

$$F_t(x) = \frac{1}{t} \int_0^t \mathbb{1}_{(l,x]}(X_u) \, du.$$

Kutoyants [13], Negri [18] and Kutoyants and Negri [15] studied this object for a certain class of stochastic differential equations. In particular, Negri [18] proved that for these particular models, $\sqrt{t}(F_t - F)$ converges weakly to a Gaussian limit, uniformly on the entire state space. We extend their results to general regular diffusions and obtain necessary and sufficient conditions in terms of the scale function and stationary distribution.

In our general setting, it follows from the classical central limit theorem that the finite-dimensional distributions of $\sqrt{t}(F_t - F)$ converge weakly to those of a centered, Gaussian random map $\mathbb{H}$ with covariance function

$$\mathrm{E}\mathbb{H}(x)\mathbb{H}(y) = 4m(I) \int_I (F(u \wedge x) - F(u)F(x))(F(u \wedge y) - F(u)F(y)) \, ds(u).$$

For uniform weak convergence we can give a necessary and sufficient integrability condition, analogous to the preceding result for the local time estimator.

By the occupation times formula (1.4)

$$F_t(x) - F(x) = \int_{(l,x]} (l_t/t - 1/m(I)) \, dm.$$

This suggests that a version of the limit process $\mathbb{H}$ is given by the process $\mathbb{H}(x) = \int_l^x \mathbb{G} \, dm$ for $\mathbb{G}$ the limit process of the diffusion local time process obtained in Theorem 2.1. In the proof of the following theorem it is seen that this integral is indeed well defined, in an $L^2$-sense, and gives a version of $\mathbb{H}$.

THEOREM 2.9. (i) *We have* $\sqrt{t}(F_t - F) \overset{P_x}{\Longrightarrow} \mathbb{H}$ *in* $\ell^\infty(J)$ *for every compact* $J \subseteq I$ *and some* $x \in I$ *if and only if* $\int_I F^2(1-F)^2 \, ds < \infty$.

(ii) *We have* $\sqrt{t}(F_t - F) \overset{P_x}{\Longrightarrow} \mathbb{H}$ *in* $\ell^\infty(I)$ *for some* $x \in I$ *if and only if* $\mathbb{H}$ *admits a version such that* $\mathbb{H}(x) \to 0$ *almost surely as* $x \downarrow l$ *or* $x \uparrow r$.

PROOF. For any $x \in I$, we have

$$E\mathbb{H}^2(x) = 4m(I)\Big((1 - F(x))^2 \int_l^x F^2 \, ds + F^2(x) \int_x^r (1-F)^2 \, ds\Big).$$

Therefore, the integrability of the function $F^2(1-F)^2$ relative to $s$ is equivalent to the existence of the covariance process of the limit process $\mathbb{H}$. It is clearly necessary for both (i) and (ii).

If $x_n$ is such that $s(x_n) = e^n$, then the integrability and monotonicity of $(1-F)^2$ near $r$ imply that $\sum_n (1-F)^2(x_n)(s(x_n) - s(x_{n-1})) < \infty$. Because $s(x_n)/s(x_{n-1}) = e$, this implies that $(1-F)^2(x_n)s(x_n) \to 0$. Again by monotonicity of $F$ we obtain that $(1-F)^2(x)s(x) \to 0$ as $x \uparrow r$. Similarly $F^2(x)s(x) \to 0$ as $x \downarrow l$.

The process $\mathbb{G}$ of Theorem 2.1 possesses continuous sample paths and hence is integrable on compacts $J \subset I$. By straightforward calculations we see that, for $a < b$ in $I$,

$$\frac{m(I)}{4} E\Big(\int_a^b \mathbb{G} \, dF\Big)^2$$

(2.2)
$$= \int_a^b \int_a^b \int (\mathbb{1}_{\{x \leq u\}} - F(u))(\mathbb{1}_{\{y \leq u\}} - F(u)) \, ds(u) \, dF(x) \, dF(y)$$

$$= (F(b) - F(a))^2 \Big(\int_l^a F^2 \, ds + \int_b^r (1-F)^2 \, ds\Big)$$

$$+ \int_a^b (F(1 - F(b)) - (1-F)F(a))^2 \, ds.$$

DONSKER THEOREMS FOR DIFFUSIONS 15

The last integral on the right-hand side is bounded above by $2(1-F(b))^2(|s(b)|+C) + 2F^2(a)(|s(a)|+C)$ for a constant $C$ [depending on $\int F^2(1-F)^2 \, ds$]. Combined with the result of the preceding paragraph and the assumed integrability of the function $F^2(1-F)^2$ it follows that $\int_a^b \mathbb{G} \, dm \to 0$ in $L^2$ as $a \to l$ and $b \to r$. Similarly, the same is true if both $a \to l$ and $b \to l$, whence the integral $\mathbb{H}(b) = \int_l^b \mathbb{G} \, dm$ is well defined in the $L^2$-sense. It can be checked that it gives a version of the limit process $\mathbb{H}$.

(i) It suffices to prove that there exists a version of the limit process $\mathbb{H}$ with sample paths that are bounded and $d_{\mathbb{H}}$-uniformly continuous on the compact $J \subset I$. In view of the preceding we have that

$$\mathrm{E}(\mathbb{H}(a) - \mathbb{H}(b))^2 = \frac{2}{\pi}(\mathrm{E}|\mathbb{H}(a) - \mathbb{H}(b)|)^2 \lesssim \left( \int_a^b \mathrm{E}|\mathbb{G}| \, dF \right)^2$$

$$\leq \left( \sup_{a<u<b} \mathrm{E}|\mathbb{G}(u)| \right)^2 (F(b) - F(a))^2$$

$$\lesssim (1 + |s(a)| \vee |s(b)|)(F(b) - F(a))^2,$$

by (2.1). It follows that for every given compact $J \subset I$ there exists a constant $C$ such that $d_{\mathbb{H}}(x,y) \leq C|F(x) - F(y)|$ for all $x, y \in J$. Since $F$ maps $J$ into the compact interval $[0,1]$, this implies that $(J, d_{\mathbb{H}})$ has finite entropy integral and hence $\mathbb{H}$ admits a version with bounded and uniformly continuous sample paths on $J$. [We can also apply Corollary 1.4 to see directly that $\sqrt{t}(F_t - F) \Rightarrow \mathbb{H}$ in $\ell^\infty(J)$ if $J$ is compact.]

(ii) Because the sample paths of the processes $\sqrt{t}(F_t - F)$ tend to zero at the boundary points of $I$, this must be true also for the limit process $\mathbb{H}$. Therefore, the existence of a version with this property is certainly necessary. We can argue the sufficiency in exactly the same manner as in the proof of Theorem 2.6. $\square$

The following corollary gives a simple sufficient condition for the sample paths of $\mathbb{H}$ to vanish at the boundary of $I$, as required in (ii) of the preceding theorem.

COROLLARY 2.10. *Suppose that $\int_I F^2(1-F)^2 \, ds < \infty$. If $(1-F)^2(x)s(x) \times \log \log s(x) \to 0$ as $x \uparrow r$ and $F^2(x)s(x) \log\log|s(x)| \to 0$ as $x \downarrow l$, then the convergence*

$$\sqrt{t}(F_t - F) \stackrel{\mathrm{P}_x}{\Longrightarrow} \mathbb{H}$$

*takes place in $\ell^\infty(I)$, for every $x \in I$. If the functions $(1-F)^2 s$ and $F^2 s$ are monotone near $r$ and $l$, respectively, then these conditions are necessary.*

PROOF. It suffices to show that the sample paths of the process $\mathbb{H}$ tend to zero at the boundary points of $I$.

Choose the sequence $x_n$ such that $s(x_n) = e^n$. Then $s(x_n)/s(x_{n-1}) = e$ and hence, for $m \leq n$, with $b_n^2 = (1-F)^2(x_n)s(x_n)$,

$$\int_{x_m}^{x_n} (1-F)^2 \, ds \lesssim \sum_{k=m}^{n} (1-F)^2(x_k)s(x_k) = \sum_{k=m}^{n} b_k^2 \lesssim \sum_{k=m}^{n-1} b_k^2.$$

From (2.2) it can be seen that a multiple of the right-hand side of this equation is a bound on $\mathrm{E}(\mathbb{H}(x_n) - \mathbb{H}(x_m))^2$.

By the bounds given in the preceding proof, for $a, b \in [x_{n-1}, x_n]$,

$$\mathrm{E}(\mathbb{H}(a) - \mathbb{H}(b))^2 \lesssim (F(b) - F(a))^2 s(x_n) =: e_n^2(a, b).$$

In particular, for $x \in [x_{n-1}, x_n]$ we have that $\mathrm{E}(\mathbb{H}(x) - \mathbb{H}(x_n))^2 \lesssim b_{n-1}^2$. It also follows that

$$N(\varepsilon, [x_{n-1}, x_n], e_n) \leq N\left(\frac{\varepsilon}{\sqrt{s(x_n)}}, [F(x_{n-1}), F(x_n)], |\cdot|\right) \lesssim \frac{b_{n-1}}{\varepsilon}.$$

Therefore, by Talagrand [28], for all $\lambda > 0$ and sufficiently large $n$ and some constant $C$,

$$(2.3) \qquad \mathrm{P}\left(\sup_{x_{n-1} \leq x \leq x_n} (\mathbb{H}(x) - \mathbb{H}(x_n)) \geq \lambda\right) \lesssim e^{-C\lambda^2/b_{n-1}^2}.$$

If $b_n^2 \log n \to 0$, then the series obtained by summing the right-hand side over $n$ is convergent for any $\lambda > 0$. In view of the definitions of $b_n$ and $x_n$ this is the case under the condition of the corollary. This implies that $\limsup_{n \to \infty} \sup_{x_{n-1} \leq x \leq x_n} (\mathbb{H}(x) - \mathbb{H}(x_n)) \leq 0$ almost surely, as $n \to \infty$. By a similar argument on the other tail we see that $\sup_{x_{n-1} \leq x \leq x_n} |\mathbb{H}(x) - \mathbb{H}(x_n)| \to 0$ almost surely.

Given a sequence of independent zero-mean Gaussian random variables $X_1, X_2, \ldots$ with $\operatorname{var} X_i = b_i^2$, let $W_n = \sum_{i=n}^{\infty} X_i$. Because $\sum_k b_k^2 < \infty$, the series $W_n$ converges in $L^2$ and hence also almost surely, by the Itô–Nisio theorem. Thus the variables $W_n$ form a well-defined Gaussian process and $W_n \to 0$ almost surely as $n \to \infty$. As noted in the preceding we have that $\mathrm{E}(\mathbb{H}(x_n) - \mathbb{H}(x_m))^2 \lesssim \sum_{k=m}^{n-1} b_k^2 = \mathrm{E}(W_n - W_m)^2$ for every $n, m \in \mathbb{N}$. This inequality remains true for $m, n \in \mathbb{N} \cup \{\infty\}$ if we set $\mathbb{H}(x_\infty) = W_\infty = 0$. Therefore, by Slepian's lemma,

$$\mathrm{E}\sup_{k \geq n} |\mathbb{H}(x_k)| \leq \mathrm{E}\sup_{\infty \geq k, l \geq n} (\mathbb{H}(x_k) - \mathbb{H}(x_l)) \leq 2\mathrm{E}\sup_{k \geq n} \mathbb{H}(x_k) \lesssim \mathrm{E}\sup_{k \geq n} W_k.$$

Because the sequence $\sup_{k \geq n} |W_k|$ converges to zero in probability as $n \to \infty$, its sequence of medians converges to zero. In view of Borell's inequality the same is then true for the sequence of means. Combined with the preceding display this shows that $\sup_{k\geq n}|\mathbb{H}(x_k)|$ converges to zero in probability, and hence $\mathbb{H}(x_n)\to 0$ almost surely.

By combining the results of the two preceding paragraphs we see that $\sup_{x\geq x_n}|\mathbb{H}(x)|\to 0$ almost surely. A similar argument applies to the limit of $\mathbb{H}$ at the left boundary of $I$. This concludes the proof of sufficiency of the condition for the Donsker property.

If the function $(1-F)^2 s$ is decreasing near $r$, then $1-F(x_n)\leq e^{(m-n)/2}(1-F(x_m))$ for $n>m$ large enough and hence $F(x_n)-F(x_m)\geq (1-F(x_m))(1-e^{-1/2})$. From (2.2) it follows that, for $n>m$ and sufficiently large $m$,

$$\mathrm{E}(\mathbb{H}(x_n)-\mathbb{H}(x_m))^2 \geq (F(x_n)-F(x_m))^2 \int_l^{x_m} F^2\, ds \gtrsim (1-F(x_m))^2 s(x_m).$$

Arguing as in the proof of Corollary 2.7 this yields the necessity of the right tail condition. The condition on the left tail can be seen to be necessary in the same way. $\square$

Because the set of indicator functions of cells in the real line is a VC-class, we can deduce the assertion of the preceding corollary also from Theorem 2.2 under the condition

$$\int_I \sqrt{|s|}\, dm < \infty.$$

For distribution functions $F$ and scale functions $s$ with regular tail behavior this condition appears to be generally stronger than the condition of the preceding corollary. For instance, if $s(x)=x$ and $1-F(x)=x^{-1/2}(\log x)^{-\alpha}$ for large $x$, then the right tail of the integral in the preceding display is finite if $\alpha>1$, whereas $(1-F)^2(x)s(x)\log\log s(x)\to 0$ as $x\to\infty$ for any $\alpha>0$. More generally, we have the following relationships between the conditions, where we state the results for the right tails only.

LEMMA 2.11. *Suppose that $\int_I \sqrt{|s|}\, dm < \infty$. Then:*

(i) $\int_I F^2(1-F)^2\, ds < \infty$.
(ii) $(1-F)^2(x)s(x)\to 0$ as $x\uparrow r$.
(iii) *If* $(1-F)^2(x)s(x)\downarrow 0$ as $x\uparrow r$, *then* $\int_x^\infty (1-F)^2\, ds\log s(x)\to 0$.
(iv) *If* $\int_x^\infty (1-F)^2\, ds\log\log s(x)\to 0$, *then* $(1-F)^2(x)s(x)\log\log s(x)\to 0$.

PROOF. By Markov's inequality we obtain, with $X_t$ a stationary diffusion, for $x$ such that $s(x)>0$,

$$1-F(x) = \mathrm{P}(\sqrt{s(X_t)}>\sqrt{s(x)}\,) \leq \frac{1}{\sqrt{s(x)}}\int_x^r \sqrt{s}\, dF.$$

In particular, the function $(1-F)\sqrt{s}$ tends to zero at the right endpoint of $I$, proving (ii). Then partial integration gives that, for $x_0$ such that $s(x_0) = 0$,

$$(2.4) \qquad \int_{x_0}^r \sqrt{s}\, dF = \frac{1}{2} \int_{x_0}^r \frac{1}{\sqrt{s}} (1-F)\, ds.$$

We conclude that finiteness of the two integrals in the display is equivalent.

(i) Because $F^2(1-F)^2 \lesssim (1-F)/\sqrt{s}$ we obtain that $\int_{x_0}^r F^2(1-F)^2\, ds < \infty$. Convergence of this integral at the left endpoint of $I$ is proved similarly.

(iii) Define $x_n$ by $s(x_n) = e^n$. Integrability of the function $(1-F)/\sqrt{s}$ at the right end of $I$ implies that $\sum_n (1-F)(x_n) e^{n/2} < \infty$. Because the sequence $(1-F)(x_n) e^{n/2}$ is decreasing by assumption, it follows that $(1-F)(x_n) e^{n/2} n \to 0$. (Indeed, if $\sum a_n < \infty$ and $a_n$ is decreasing, then $\infty > \sum_k \sum_{2^{k-1} \le n < 2^k} a_n \ge \sum_k 2^{k-1} a_{2^{k-1}}$, so that $2^k a_{2^k} \to 0$ as $k \to \infty$. It follows that $\sup_{2^{k-1} \le n \le 2^k} n a_n \lesssim 2^{k-1} a_{2^{k-1}} \to 0$, so $n a_n \to 0$.) Hence,

$$\sum_{n \ge n_0} (1-F)^2(x_n) s(x_n) \le (1-F)(x_{n_0}) e^{n_0/2} \sum_{n \ge n_0} (1-F)(x_n) e^{n/2} = O(1/n_0),$$

as $n_0 \to \infty$. This implies that $\log s(x_{n_0}) \int_{x_{n_0}}^r (1-F)^2\, ds \to 0$.

(iv) With $x_n$ as before, we have $(1-F)^2(x_{n_0}) s(x_{n_0}) \le \sum_{n \ge n_0} (1-F)^2(x_n) \times s(x_n)$, which is bounded above by a multiple of $\int_{x_{n_0-1}}^r (1-F)^2\, ds$. □

On the other hand, it is not true in general that the condition $\int_I \sqrt{|s|}\, dm < \infty$ is stronger than the condition of Corollary 2.10 and hence the latter condition is not necessary in general. This is also clear from the proof, which is based on the assumption that the right-hand side of (2.3) yields a convergent series. Without some regularity on the sequence $b_n^2$, this does not reduce to the simple condition as stated.

EXAMPLE 2.12. Define a sequence $x_n$ by $\log \log \log \log \log s(x_n) = n$ (where we use the logarithm at base 2), and define

$$1 - F(x) = \frac{1}{\sqrt{s(x_n)} \sqrt{\log \log \log \log s(x_n)}}, \qquad x_{n-1} < x \le x_n.$$

Then $\int_I \sqrt{|s|}\, dm < \infty$, but $(1-F)^2(x_n)\, s(x_n) = (\log \log \log \log s(x_n))^{-1}$.

Because this distribution function $F$ possesses flat parts, it cannot appear as the stationary distribution of a regular diffusion. However, by moving a tiny fraction of the total mass, we can construct a distribution with full support without destroying the preceding properties.

**3. Continuous martingales and majorizing measures.** Let $(\Omega, \mathcal{F}, \{\mathcal{F}_t\}, \mathrm{P})$ be a filtered probability space. On this stochastic basis, suppose that we have a collection $M = \{M^\theta : \theta \in \Theta\}$ of continuous local martingales $M^\theta = (M^\theta_t : t \geq 0)$, indexed by a countable pseudo-metric space $(\Theta, d)$. The *quadratic d-modulus of continuity* $\|M\|_d$ of the collection $M$ is the stochastic process defined by

$$\|M\|_{d,t} = \sup_{\theta, \psi \,:\, d(\theta,\psi) > 0} \frac{\sqrt{\langle M^\theta - M^\psi \rangle_t}}{d(\theta, \psi)}.$$

Here $\langle N \rangle$ denotes the quadratic variation process of the continuous local martingale $N$.

The quadratic modulus was introduced explicitly by Nishiyama [21, 22] and appeared already implicitly in the papers Bae and Levental [2] and Nishiyama [20]. The relevance of the quadratic modulus stems from the fact that for every time $t \geq 0$ and every constant $K > 0$, the random map $\theta \mapsto M^\theta_t \mathbb{1}_{\{\|M\|_{d,t} \leq K\}}$ is sub-Gaussian with respect to the pseudo-metric $Kd$. Indeed, the Bernstein inequality for continuous local martingales (see, e.g., [26]) implies that

$$\begin{aligned}
&\mathrm{P}(|M^\theta_t \mathbb{1}_{\{\|M\|_{d,t} \leq K\}} - M^\psi_t \mathbb{1}_{\{\|M\|_{d,t} \leq K\}}| \geq x) \\
&\quad \leq \mathrm{P}(|M^\theta_t - M^\psi_t| \geq x, \|M\|_{d,t} \leq K) \\
&\quad \leq \mathrm{P}(|M^\theta_t - M^\psi_t| \geq x, \langle M^\theta - M^\psi \rangle_t \leq K^2 d^2(\theta, \psi)) \\
&\quad \leq 2 e^{-(1/2) x^2 / (K^2 d^2(\theta, \psi))}.
\end{aligned}$$

For random maps whose increments are controlled in this manner, the theory of majorizing measures gives sharp bounds for the modulus of continuity. As before, we denote by $N(\eta, \Theta, d)$ the minimal number of balls of $d$-radius $\eta$ that are needed to cover $\Theta$. The symbol $\lesssim$ between two expressions means that the left-hand side is less than a universal positive constant times the right-hand side.

LEMMA 3.1. *For all $\delta, x, \eta > 0$, $K \geq 1$, every Borel probability measure $\nu$ on $(\Theta, d)$, and every bounded stopping time $\tau$,*

$$\mathrm{P}\left( \sup_{t \leq \tau} \sup_{d(\theta, \psi) < \delta} |M^\theta_t - M^\psi_t| \geq x; \|M\|_{d,\tau} \leq K \right)$$

$$\lesssim \frac{K}{x} \left( \sup_\theta \int_0^\eta \sqrt{\log \frac{1}{\nu(B_d(\theta, \varepsilon))}} \, d\varepsilon + \delta \sqrt{N(\eta, \Theta, d)} \right).$$

PROOF. We may of course assume that the right-hand side of the inequality in the statement of the lemma is finite. Introduce the stopping time $\tau_K = \inf\{t : \|M\|_{d,t} > K\}$, so that the probability in the statement of the lemma is bounded by $\mathrm{P}(\sup_{t \leq \tau} X_t^\delta \geq x)$, where

$$(3.1) \qquad X_t^\delta = \sup_{d(\theta,\psi) < \delta} |M_{\tau_K \wedge t}^\theta - M_{\tau_K \wedge t}^\psi|.$$

By Bernstein's exponential inequality for continuous martingales we have for all $a \geq 0$ and every finite stopping time $\sigma$

$$\begin{aligned}
& \mathrm{P}(|M_{\tau_K \wedge \sigma}^\theta - M_{\tau_K \wedge \sigma}^\psi| > a) \\
& = \mathrm{P}(|M_{\tau_K \wedge \sigma}^\theta - M_{\tau_K \wedge \sigma}^\psi| > a; \langle M^\theta - M^\psi \rangle_{\tau_K \wedge \sigma} \leq K^2 d^2(\theta,\psi)) \\
& \leq 2 e^{-(1/2)a^2/(K^2 d^2(\theta,\psi))}.
\end{aligned}$$

Hence, the random map $\theta \mapsto M_{\tau_K \wedge \sigma}^\theta$ is sub-Gaussian with respect to the pseudo-metric $Kd$. By formula (11.15) on page 317 of [16] this implies that for all $\delta, \eta > 0$

$$(3.2) \qquad \mathrm{E} X_\sigma^\delta \lesssim K \left( \sup_\theta \int_0^\eta \sqrt{\log \frac{1}{\nu(B_d(\theta,\varepsilon))}} \, d\varepsilon + \delta \sqrt{N(\eta, \Theta, d)} \right),$$

where $B_d(\xi, \varepsilon)$ is the ball around $\xi$ with $d$-radius $\varepsilon$. In particular, we see that $\mathrm{E} X_t^\delta < \infty$ for every $t \geq 0$. Also, for any pair $(\theta, \psi)$ and for every finite stopping time $\sigma$, by the Davis–Gundy inequality,

$$\mathrm{E}(M_{\tau_K \wedge \sigma}^\theta - M_{\tau_K \wedge \sigma}^\psi)^2 \leq \mathrm{E} \langle M^\theta - M^\psi \rangle_{\tau_K \wedge \sigma} \leq K^2 d^2(\theta, \psi).$$

Thus, the collection $\{M_{\tau_K \wedge \sigma}^\theta - M_{\tau_K \wedge \sigma}^\psi : \sigma \text{ is a finite stopping time}\}$ is bounded in $L^2$ and therefore uniformly integrable. This implies that the stopped local martingale $M_{\tau_K \wedge t}^\theta - M_{\tau_K \wedge t}^\psi$ is of class (D), which means that it is in fact a uniformly integrable martingale (see, e.g., pages 11–12 of [10]). It is then easy to see that the process $X^\delta$ defined by (3.1) is a submartingale. Hence, by the submartingale inequality, $\mathrm{P}(\sup_{t \leq \tau} X_t^\delta \geq x) \leq \mathrm{E} X_\tau^\delta / x$. In combination with (3.2) this yields the statement of the lemma. □

With the help of this lemma we can prove results concerning the regularity and asymptotic tightness of collections of continuous local martingales under majorizing measure conditions. The key condition is the existence of a pseudo-metric $d$ on $\Theta$ for which the modulus is finite or bounded in probability and for which there exists a probability measure $\nu$ such that the integral on the right-hand side in the preceding lemma is continuous at zero. The latter is the continuous majorizing measure condition:

$$(3.3) \qquad \lim_{\eta \downarrow 0} \sup_\theta \int_0^\eta \sqrt{\log \frac{1}{\nu(B_d(\theta,\varepsilon))}} \, d\varepsilon = 0.$$

The first theorem deals with regularity of a given collection of local martingales $M$.

THEOREM 3.2. *Suppose there exists a Borel probability measure $\nu$ on $(\Theta, d)$ such that (3.3) holds for a pseudo-metric $d$ on $\Theta$ for which $\|M\|_{d,\tau} < \infty$ almost surely. Then the random map $\theta \mapsto M_\tau^\theta$ is almost surely bounded and uniformly $d$-continuous on $\Theta$.*

PROOF. By Lemma 3.1, there exists for every $n \in \mathbb{N}$ a positive number $\delta_n$ such that for every $K, x > 0$,

$$P\left(\sup_{d(\theta,\psi)<\delta_n} |M_\tau^\theta - M_\tau^\psi| \geq x; \|M\|_{d,\tau} \leq K\right) \lesssim \frac{K}{4^n x}.$$

For every $n$, define the event

$$A_n = \left\{\sup_{d(\theta,\psi)<\delta_n} |M_\tau^\theta - M_\tau^\psi| > \frac{1}{2^n}\right\}.$$

Then for every $K > 0$ we have $\sum P(A_n; \|M\|_{d,\tau} \leq K) \lesssim K \sum 2^{-n} < \infty$. So by the Borel–Cantelli lemma, $P(A_n \text{ infinitely often}; \|M\|_{d,\tau} \leq K) = 0$. Since $\|M\|_{d,\tau}$ is almost surely finite by assumption, it follows that

$$P(A_n \text{ infinitely often}) = P(A_n \text{ infinitely often}; \|M\|_{d,\tau} < \infty)$$

$$\leq \sum_K P(A_n \text{ infinitely often}; \|M\|_{d,\tau} \leq K) = 0.$$

So we almost surely have that $\sup_{d(\theta,\psi)<\delta_n} |M_\tau^\theta - M_\tau^\psi| \leq 2^{-n}$ for all $n$ large enough, which implies that the random map $\theta \mapsto M_\tau^\theta$ is uniformly continuous. Recall that under the majorizing measure condition, the pseudo-metric space $(\Theta, d)$ is totally bounded (see, e.g., the proof of Lemma A.2.19 of [29]). It follows that $\theta \mapsto M_\tau^\theta$ is bounded with probability 1. □

Suppose now that for each $n \in \mathbb{N}$, we have a collection $M^n = \{M^{n,\theta} : \theta \in \Theta\}$ of continuous local martingales and a finite stopping time $\tau_n$ on a stochastic basis $(\Omega^n, \mathcal{F}^n, \{\mathcal{F}_t^n\}, P^n)$. For each $n$ the local martingales $M^{n,\theta}$ are indexed by a parameter $\theta$ belonging to a fixed pseudo-metric space $\Theta$. Recall that a sequence $X_n$ of $\ell^\infty(\Theta)$-valued random elements is called *asymptotically $d$-equicontinuous in probability* if for all $\varepsilon, \eta > 0$ there exists a $\delta > 0$ such that

$$\limsup_{n \to \infty} P\left(\sup_{d(\theta,\psi)\leq\delta} |X_n(\theta) - X_n(\psi)| > \varepsilon\right) \leq \eta.$$

Weak convergence in $\ell^\infty(\Theta)$ to a tight limit is equivalent to finite-dimensional convergence and equicontinuity with respect to a semimetric $d$ such that $(\Theta, d)$ is totally bounded (see, e.g., [29], Theorem 1.5.7). For the random maps $\theta \mapsto M^{n,\theta}_{\tau_n}$, finite-dimensional weak convergence will typically follow from a classical martingale central limit theorem (cf. [10]). Using Lemma 3.1, it is straightforward to give sufficient conditions for asymptotic equicontinuity in terms of the quadratic modulus and majorizing measures. The next theorem extends Theorem 3.2.4 of [22], which gives sufficient conditions for asymptotic equicontinuity in terms of metric entropy.

THEOREM 3.3. *Suppose there exists a Borel probability measure $\nu$ on $(\Theta, d)$ such that (3.3) holds for a pseudo-metric $d$ on $\Theta$ for which $\|M^n\|_{d,\tau_n} = O_P(1)$. Then $(\Theta, d)$ is totally bounded and the sequence of random maps $\theta \mapsto M^{n,\theta}_{\tau_n}$ in $\ell^\infty(\Theta)$ is asymptotically $d$-equicontinuous in probability.*

PROOF. The total boundedness of $(\Theta, d)$ is a direct consequence of the existence of a majorizing measure. See, for example, the proof of Lemma A.2.19 of [29].

Let the random map $X_n$ on $\Theta$ be defined by $X_n(\theta) = M^{n,\theta}_{\tau_n}$. Then for every $K > 0$

$$P\bigg(\sup_{d(\theta,\psi)\leq \delta} |X_n(\theta) - X_n(\psi)| > \varepsilon\bigg)$$
$$\leq P\bigg(\sup_{d(\theta,\psi)\leq \delta} |X_n(\theta) - X_n(\psi)| > \varepsilon; \|M^n\|_{d,\tau_n} \leq K\bigg) + P(\|M^n\|_{d,\tau_n} > K).$$

Now if $\eta > 0$ is given, we can first choose $K$ large enough to ensure that $\limsup P(\|M^n\|_{d,\tau_n} > K) < \eta/2$. Lemma 3.1 implies that for this fixed $K$, we can choose a $\delta > 0$ such that the first term on the right-hand side is less than $\eta/2$. □

The preceding theorems do not use the full power of Lemma 3.1, because they use the control in $\theta$ of the local martingales $t \mapsto M^{n,\theta}_t$, but not the control in the time parameter $t$. In the following theorem we use the lemma to establish a majorizing measure condition for the asymptotic tightness in $\ell^\infty([0,T] \times \Theta)$ of random maps of the form $(t, \theta) \mapsto M^{n,\theta}_t$, for fixed $T \in (0, \infty)$.

We make the same assumptions as in the preceding theorem, and in addition assume that for every fixed $\theta \in \Theta$ the sequence of processes $(M^{n,\theta}_t : 0 \leq t \leq T)$ is asymptotically equicontinuous in probability relative to the Euclidean metric on $[0,T]$. By the martingale central limit theorem, this is certainly true if the sequence of quadratic variation processes $\langle M^{n,\theta} \rangle$ converges pointwise in probability to a continuous function (which is then the quadratic variation process of the Gaussian limit process).

THEOREM 3.4. *Suppose there exists a Borel probability measure $\nu$ on $(\Theta, d)$ such that (3.3) holds for a pseudo-metric $d$ on $\Theta$ for which $\|M^n\|_{d,\tau_n} = O_P(1)$. Furthermore, assume that, for every fixed $\theta \in \Theta$, the sequence of processes $(M_t^{n,\theta} : 0 \leq t \leq T)$ is asymptotically equicontinuous in probability relative to the Euclidean metric. Then the sequence of random maps $M^n$ is asymptotically tight in the space $\ell^\infty([0,T] \times \Theta)$.*

PROOF. By the majorizing measure condition (3.3) the set $\Theta$ is totally bounded under $d$. If $\theta_1, \ldots, \theta_m$ is a $\delta$-net over $\Theta$ and $s, t \in [0, T]$, then for all $i$

$$|M_s^{n,\theta} - M_t^{n,\theta}| \leq |M_s^{n,\theta_i} - M_t^{n,\theta_i}| + 2 \sup_{0 \leq t \leq T} |M_t^{n,\theta} - M_t^{n,\theta_i}|.$$

Hence

$$\sup_{|s-t|<\gamma} \sup_{d(\theta,\psi)\leq\delta} |M_s^{n,\theta} - M_t^{n,\psi}|$$

$$\leq \sup_{|s-t|<\gamma} \sup_{d(\theta,\psi)\leq\delta} (|M_s^{n,\theta} - M_t^{n,\theta}| + |M_t^{n,\theta} - M_t^{n,\psi}|)$$

$$\leq \max_i \sup_{|s-t|<\gamma} |M_s^{n,\theta_i} - M_t^{n,\theta_i}| + 3 \sup_{0 \leq t \leq T} \sup_{d(\theta,\psi)\leq\delta} |M_t^{n,\theta} - M_t^{n,\psi}|.$$

Fix $\varepsilon, \eta > 0$. Extending the argument in the proof of Theorem 3.3, we can show that there exists $\delta > 0$ such that

$$(3.4) \qquad \limsup_{n \to \infty} P\left( \sup_{0 \leq t \leq T} \sup_{d(\theta,\psi)\leq\delta} |M_t^{n,\theta} - M_t^{n,\psi}| > \varepsilon \right) < \eta.$$

For this $\delta = \delta(\varepsilon, \eta)$ there exists a finite $\delta$-net $\theta_1, \ldots, \theta_m$ over $\Theta$ (where $m$ depends on $\delta$). By the assumption of asymptotic equicontinuity of the processes $t \mapsto M_t^{n,\theta}$, there exists $\gamma = \gamma(\eta, m, \theta_1, \ldots, \theta_m)$ such that

$$\limsup_{n\to\infty} P\left( \sup_{|s-t|<\gamma} |M_s^{n,\theta_i} - M_t^{n,\theta_i}| > \varepsilon \right) < \frac{\eta}{m}, \qquad i = 1, \ldots, m.$$

Combining the preceding displays we see that

$$\limsup_{n\to\infty} P\left( \sup_{|s-t|<\gamma} \sup_{d(\theta,\psi)\leq\delta} |M_s^{n,\theta} - M_t^{n,\psi}| > 4\varepsilon \right) \leq \sum_{i=1}^m \frac{\eta}{m} + \eta \leq 2\eta.$$

Thus for the given pair $(\varepsilon, \eta)$ we have found a pair $(\gamma, \delta)$ of positive numbers such that this holds. Because the probability on the left-hand side is increasing in $\gamma$ and $\delta$, the bound remains true if we replace $\gamma$ or $\delta$ by the smaller of the two. This implies that the sequence of processes $M^n$ is asymptotically equicontinuous in probability relative to the product of the Euclidean metric on $[0, T]$ and the pseudo-metric $d$ on $\Theta$, and hence it is asymptotically tight ([29], Theorem 1.5.7). □

In the preceding theorem we can also use an arbitrary pseudo-metric for which the interval $[0, T]$ is totally bounded (and this could be permitted to depend on $\theta$), rather than the Euclidean metric. However, because the local martingales $t \mapsto M_t^{n,\theta}$ are continuous relative to the Euclidean metric by assumption, this apparent generalization would not make the theorem more general: the necessary continuity of the limit points $t \mapsto M_t^\theta$ would imply that the equicontinuity necessarily also holds relative to the Euclidean pseudo-metric. For simplicity of the statement we have used the Euclidean metric throughout.

**4. A limit theorem for diffusion local time.** In this section we collect some classical and some less well-known facts about diffusion local time. We shall need these in the proof of Theorem 1.2. As in the Introduction, let $X$ be the regular diffusion on the open interval $I$. A central result in the theory of one-dimensional diffusions is that diffusions in natural scale are in fact time-changed Brownian motions; see, for instance, [25] or [11]. In our setting, we have that under $P_x$, it holds that $s(X_t) = W_{\tau_t}$, where $W$ is a Brownian motion that starts in $s(x)$, and $\tau_t$ is the right-continuous inverse of the process $A$ defined by

$$A_t = \int_I L_t^W(s(y)) m(dy).$$

Here $L^W = (L_t^W(y) : t \geq 0, y \in \mathbb{R})$ is the local time of $W$. It follows from this relation that the local time $l_t(y)$ of $X$ with respect to the speed measure $m$ satisfies $l_t(y) = L_{\tau_t}^W(s(y))$.

This time-change representation of diffusion local time shows that with probability 1, the random function $y \mapsto l_t(y)$ can be chosen continuous and has compact support. In particular, it holds that $\|l_t\|_\infty = \sup_{y \in I} l_t(y) < \infty$ almost surely. In [30] it is shown that in fact, $\|l_t\|_\infty = O_P(t)$ as $t \to \infty$. For the sake of easy reference, we include a proof of this fact. We need the following lemma.

LEMMA 4.1. *For every $x \in I$ we have, for $Z$ standard normally distributed and $t \to \infty$,*

$$\frac{\tau_t}{t^2} \stackrel{P_x}{\Longrightarrow} \frac{1}{m^2(I) Z^2}.$$

PROOF. The process $A$ defined above is a continuous additive functional of $W$, and since $m$ is finite, it is integrable. By Proposition (2.2) in Chapter XIII of [24], it follows that

$$\text{(4.1)} \qquad \frac{1}{\sqrt{t}} A_t \xrightarrow{\mathrm{P}_x} m(I) L_1^B(0),$$

where $L^B$ is the local time of a Brownian motion $B$ that starts in 0. The process $\tau$ is the right-continuous inverse of $A$, so for every $t, T \geq 0$ it holds that $\tau_t < T$ if and only if $A_T > t$. By (4.1), it follows that, for every $z \geq 0$,

$$\begin{aligned}
\mathrm{P}_x\left(\frac{\tau_t}{t^2} < z\right) &= \mathrm{P}_x(A_{t^2 z} > t) = \mathrm{P}_x\left(\frac{1}{t\sqrt{z}} A_{t^2 z} > \frac{1}{\sqrt{z}}\right) \\
&\to \mathrm{P}_x\left(m(I) L_1^B(0) > \frac{1}{\sqrt{z}}\right) \\
&= \mathrm{P}_x\left(\frac{1}{m^2(I)(L_1^B(0))^2} < z\right).
\end{aligned}$$

To complete the proof we use the well-known fact that $(L_1^B(0))^2$ has a $\chi_1^2$-distribution (see [12], Theorem 3.6.17 and Problem 2.8.2). □

THEOREM 4.2. *For every $x \in I$ we have $\|l_t\|_\infty = O_{\mathrm{P}_x}(t)$ as $t \to \infty$.*

PROOF. Let us write $\alpha_t = t^{-1} \|l_t\|_\infty$. We have to prove that $\alpha_t$ is asymptotically tight for $t \to \infty$. By the time-change relation, we have for all $a, b > 0$

$$\begin{aligned}
\mathrm{P}_x(\alpha_t > a) &= \mathrm{P}_x\left(\sup_{z \in s(I)} \frac{1}{t} L_{t^2(\tau_t/t^2)}^W(z) > a\right) \\
&\leq \mathrm{P}_x\left(\sup_{z \in \mathbb{R}, u \leq b} \frac{1}{t} L_{t^2 u}^W(z) > a\right) + \mathrm{P}_x\left(\frac{\tau_t}{t^2} > b\right).
\end{aligned}$$

By the scaling property of Brownian local time (see Exercise (2.11) in Chapter VI of [24] and note that $W$ is a Brownian motion starting at $s(x)$) it holds under $\mathrm{P}_x$ that

$$\sup_{z \in \mathbb{R}, u \leq b} \frac{1}{t} L_{t^2 u}^W(z) \stackrel{\mathrm{d}}{=} \sup_{z \in \mathbb{R}, u \leq b} L_u^B\left(\frac{z - s(x)}{t}\right) = \sup_{z \in \mathbb{R}} L_b^B(z),$$

where $L^B$ is the local time of a standard Brownian motion $B$ (starting in 0). So we find that for all $a, b > 0$

$$\text{(4.2)} \qquad \mathrm{P}_x(\alpha_t > a) \leq \mathrm{P}_x\left(\sup_{z \in \mathbb{R}} L_b^B(z) > a\right) + \mathrm{P}_x\left(\frac{\tau_t}{t^2} > b\right).$$

The proof is finished upon noting that $z \mapsto L_b^B(z)$ is bounded (because continuous with compact support), almost surely and $\tau_t/t^2$ is asymptotically tight. □

## 5. Proof of Theorem 1.2.

5.1. *Reduction to the natural scale case.* Let us first show that it suffices to prove the theorem for diffusions $X$ that are in natural scale (i.e., for which the identity function is a scale function). The diffusion $Y = s(X)$ is in natural scale (see, e.g., Theorem V.46.12 of [25]), and we have the relations

$$m = m^Y \circ s, \qquad l_t = l_t^Y \circ s, \qquad F = F^Y \circ s,$$

between the local time $l^Y$, speed measure $m^Y$ and stationary distribution $F^Y$ of $Y$, and the local time $l$, speed measure $m$ and stationary distribution $F$ of $X$. Moreover,

$$\mathrm{E}\mathbb{H}^Y(\lambda \circ s^{-1})\mathbb{H}^Y(\nu \circ s^{-1}) = \mathrm{E}\mathbb{H}\lambda\mathbb{H}\nu.$$

It follows that the class $\Lambda$ is Donsker for $X$ if and only if the class $\Lambda \circ s^{-1} = \{\lambda \circ s^{-1} : \lambda \in \Lambda\}$ is Donsker for $Y$. So if we have proved the theorem for diffusions in natural scale, we can apply it to the diffusion $Y = s(X)$ and the class $\Lambda \circ s^{-1}$ to prove it for a diffusion $X$ that is not in natural scale.

In the remainder of the proof we therefore assume that $X$ is in natural scale. The process $X$ is then an ergodic diffusion in natural scale on the open interval $I$. Therefore, we must have $I = \mathbb{R}$ (see, e.g., Theorem 20.15 of [11]). Moreover, the fact that the state space is open implies that $X$ is a local martingale (cf., e.g., [25], Corollary V.46.15). We also note that for diffusions in natural scale on an open interval, the diffusion local time $l_t(x)$ with respect to the speed measure coincides with the semimartingale local time of $X$ (see [25], Section V.49).

5.2. *Asymptotic equivalence with uniform weak convergence of continuous local martingales.* In this section we show that the weak convergence of the empirical process $\mathbb{H}_t$ is equivalent to the weak convergence of a normalized $\ell^\infty(\Lambda)$-valued continuous local martingale. Since $X$ is now in natural scale, we have $I = \mathbb{R}$. For every $x \in \mathbb{R}$, define the functions $\pi_x$ and $\Pi_x$ on $\mathbb{R}$ by $\pi_x = 2(\mathbb{1}_{[x,\infty)} - F)$ and

$$\Pi_x(y) = \int_{y_0}^{y} \pi_x(u)\, du,$$

where $y_0$ is an arbitrary, but fixed point in $\mathbb{R}$. The function $\pi_x$ is the difference of two increasing functions, and hence $\Pi_x$ is the difference of two convex functions. Moreover, we have the relation $\pi_x(b) - \pi_x(a) = \nu(a, b]$ for all $a \leq b$, where $\nu$ is the signed measure $\nu = 2(\delta_x - \mu)$ on $\mathbb{R}$, and $\delta_x$ denotes the Dirac measure concentrated at $x$. So by the generalized Itô formula (see, e.g., [24], Theorem VI.1.5, or [25], 45.1)

$$\Pi_x(X_t) - \Pi_x(X_0) = \int_0^t \pi_x(X_u)\, dX_u + \tfrac{1}{2}\int_{\mathbb{R}} l_t(y)\nu(dy).$$

It follows from the definition of $\nu$ and the occupation times formula (1.4) that

$$\frac{1}{2}\int_{\mathbb{R}} l_t(y)\nu(dy) = l_t(x) - \int_{\mathbb{R}} l_t(y)\mu(dy) = l_t(x) - \frac{1}{m(I)}t,$$

so that, under $P_z$,

$$\frac{1}{t}l_t(x) - \frac{1}{m(\mathbb{R})} = \frac{1}{t}(\Pi_x(X_t) - \Pi_x(z)) - \frac{1}{t}\int_0^t \pi_x(X_u)\,dX_u.$$

If we integrate this identity with respect to $\lambda(dx)$ and use the stochastic Fubini theorem (see [23], Theorem IV.45), we see that the empirical process $\mathbb{H}_t$ can be decomposed as

$$\mathbb{H}_t\lambda = R_{z,t}(\lambda) - \frac{1}{\sqrt{t}}M_t^\lambda \tag{5.1}$$

under $P_z$, where $M^\lambda$ is the continuous local martingale defined by

$$M_t^\lambda = 2\int_0^t h_\lambda(X_u)\,dX_u \quad \text{with } h_\lambda(x) = \lambda(l,x] - F(x)\lambda(I), \tag{5.2}$$

and $R_{z,t}(\lambda) = t^{-1/2}\int_{\mathbb{R}}(\Pi_x(X_t) - \Pi_x(z))\lambda(dx)$. The next step is to show that the $R_{z,t}$-term vanishes uniformly in $\lambda$, so that we only have to deal with the martingale part of $\mathbb{H}_t$. The functions $\pi_x$ are bounded in absolute value by 2, so we have the pointwise inequality $|\Pi_x| \leq \Pi$ for every $x$, where $\Pi$ is a function that does not depend on $x$. It follows that

$$\sup_{x\in\mathbb{R}}\frac{1}{\sqrt{t}}|\Pi_x(X_t) - \Pi_x(z)| \leq \frac{1}{\sqrt{t}}(\Pi(X_t) + \Pi(z)).$$

Consequently, we have

$$\sup_{\lambda\in\Lambda}|R_{z,t}(\lambda)| \leq \sup_{\lambda\in\Lambda}\|\lambda\|\frac{1}{\sqrt{t}}(\Pi(X_t) + \Pi(z)).$$

The right-hand side converges to 0 in probability, since the law of $X_t$ converges in total variation distance to the stationary measure $\mu$ as $t \to \infty$, whatever the initial law (see, e.g., [25], Section 54.5).

5.3. *Continuity of the empirical process.* In this section we prove that the empirical process is $d_{\mathbb{H}}$-continuous. Observe that we have $\mathbb{H}_t\lambda = \int_I \phi_t(x)\lambda(dx)$, with

$$\phi_t(x) = \sqrt{t}\left(\frac{1}{t}l_t(x) - \frac{1}{m(I)}\right).$$

Since the random function $x \mapsto l_t(x)$ is almost surely continuous and has compact support, the random function $\phi_t$ is bounded and continuous with probability 1. Note also that by the occupation times formula (1.4), $\phi_t$ satisfies
$$\int_I \phi_t(x) m(dx) = \frac{1}{\sqrt{t}} \int_0^t 1\, du - \sqrt{t} = 0.$$
Hence, the $d_{\mathbb{H}}$-continuity of the empirical process is a consequence of the following lemma.

LEMMA 5.1. *Let $\psi$ be a continuous function on $I$ with compact support and $\phi = \psi - \int \psi\, d\mu$. Then if $\lambda_n$ and $\lambda$ are signed measures with $\|\lambda_n\|$ bounded and $d_{\mathbb{H}}(\lambda_n, \lambda) \to 0$, it holds that*
$$\int_I \phi(x) \lambda_n(dx) \to \int_I \phi(x) \lambda(dx)$$
*as $n \to \infty$.*

PROOF. Define $\alpha_n = \lambda_n(I) - \lambda(I)$. Since the total variation of the signed measures $\lambda_n$ is uniformly bounded by assumption, the sequence $\alpha_n$ is bounded. Hence, it has a converging subsequence, say $\alpha_{n'} \to \alpha$. Observe that the convergence $d_{\mathbb{H}}(\lambda_n, \lambda) \to 0$ implies that $n'$ has a further subsequence $n''$ such that
$$\lambda_{n''}(l, x] - F(x)\lambda_{n''}(I) \to \lambda(l, x] - F(x)\lambda(I)$$
for almost every $x \in I$ [recall that $X$ is now in natural scale, so that $s(x) = x$]. So there exists a dense set $D \subseteq I$ such that for all $x \in D$
$$\lambda_{n''}(l, x] = \lambda_{n''}(l, x] - F(x)\lambda_{n''}(I) + F(x)\lambda_{n''}(I)$$
$$\to \lambda(l, x] - F(x)\lambda(I) + F(x)(\alpha + \lambda(I)) = \lambda(l, x] + \alpha F(x).$$
Since $\psi$ is compactly supported and continuous, we can approximate it uniformly by functions $\psi_m$ of the form $\psi_m = \sum_{i=1}^m \mathbb{1}_{(a_{m,i}, a_{m,i+1}]} \psi(a_{i,m})$ for $a_{m,0} < a_{m,1} < \cdots < a_{m,m}$ finite partitions of the support of $\psi$. The preceding display shows that for every fixed $m$, as $n \to \infty$,
$$\int \psi_m(x) \lambda_{n''}(dx) \to \int \psi_m(x) \lambda(dx) + \alpha \int \psi_m(x) \mu(dx).$$
By the uniform approximation as $m \to \infty$, this is then also true for $\psi$ in the place of $\psi_m$. Consequently,
$$\int_I \phi(x) \lambda_{n''}(dx) = \int_I \psi(x) \lambda_{n''}(dx) - \lambda_{n''}(I) \int_I \psi\, d\mu$$
$$\to \int_I \psi(x) \lambda(dx) + \alpha \int \psi(x) \mu(dx) - (\alpha + \lambda(I)) \int \psi(x) \mu(dx)$$
$$= \int_I \phi(x) \lambda(dx).$$
This completes the proof. □

DONSKER THEOREMS FOR DIFFUSIONS 295.4. *End of the proof.* We can now finish the proof of Theorem 1.2. That the existence of a bounded and $d_{\mathbb{H}}$-uniformly continuous version of the Gaussian limit $\mathbb{H}$ is necessary for $\Lambda$ to be Donsker follows from the general theory of weak convergence to Gaussian processes.

For the proof of the converse, suppose that such a bounded and uniformly continuous version $\mathbb{H}$ exists. By Theorem 12.9 of [16], this implies that there exists a Borel measure $\nu$ on $(\Lambda, d_{\mathbb{H}})$ such that

$$(5.3) \qquad \limsup_{\eta \downarrow 0} \sup_{\lambda} \int_0^{\eta} \sqrt{\log \frac{1}{\nu(B_{d_{\mathbb{H}}}(\lambda, \varepsilon))}} \, d\varepsilon = 0.$$

It follows that $(\Lambda, d_{\mathbb{H}})$ is totally bounded (see, e.g., the proof of Lemma A.2.19 of [29]), and therefore separable. Hence, since the empirical process $\mathbb{H}_t$ is $d_{\mathbb{H}}$-continuous, we may assume that $\Lambda$ is countable.

By the considerations in Section 5.2, it suffices to show that the random maps $\mathbb{H}'_t$ defined by $\mathbb{H}'_t \lambda = t^{-1/2} M_t^{\lambda}$ converge weakly to $\mathbb{H}$, where $M^{\lambda}$ is the local martingale given by (5.2). By the ergodic theorem ([9], Section 6.8), we have

$$\frac{1}{t} \langle M^{\lambda}, M^{\nu} \rangle_t = \frac{4}{t} \int_0^t h_{\lambda}(X_u) h_{\nu}(X_u) \, d\langle X \rangle_u = \frac{4}{t} \int_{\mathbb{R}} h_{\lambda}(x) h_{\nu}(x) l_t(x) \, dx$$

$$\stackrel{\text{a.s.}}{\to} \frac{4}{m(\mathbb{R})} \int_{\mathbb{R}} h_{\lambda}(x) h_{\nu}(x) \, dx = \mathbb{E} \mathbb{H} \lambda \mathbb{H} \nu$$

as $t \to \infty$. So by the martingale central limit theorem, the finite-dimensional distributions of $\mathbb{H}'_t$ converge weakly to those of $\mathbb{H}$. Now pick an arbitrary sequence $a_n \to \infty$ and apply Theorem 3.3 to the local martingales $M^{n,\lambda}$ defined by

$$M_t^{n,\lambda} = \frac{1}{\sqrt{a_n}} M_{a_n t}^{\lambda},$$

and the stopping time $\tau_n$ equal to 1. Then in view of (5.3), all that remains to be shown is that $\|M^n\|_{d_{\mathbb{H}},1} = O_P(1)$. Since

$$\frac{1}{a_n} \langle M^{\lambda} - M^{\nu} \rangle_{a_n} = \frac{4}{a_n} \int_{\mathbb{R}} (h_{\lambda}(x) - h_{\nu}(x))^2 l_{a_n}(x) \, dx \lesssim \frac{1}{a_n} \|l_{a_n}\|_{\infty} d_{\mathbb{H}}^2(\lambda, \nu),$$

we have

$$\|M^n\|_{d_{\mathbb{H}},1}^2 = \sup_{d_{\mathbb{H}}(\lambda,\nu)>0} \frac{(1/a_n) \langle M^{\lambda} - M^{\nu} \rangle_{a_n}}{d_{\mathbb{H}}^2(\lambda,\nu)} \lesssim \frac{1}{a_n} \|l_{a_n}\|_{\infty}.$$

By Theorem 4.2 it holds that $\|l_{a_n}\|_{\infty} = O_P(a_n)$, so indeed, $\|M^n\|_{d_{\mathbb{H}},1} = O_P(1)$.

DEPARTMENT OF MATHEMATICS
FACULTY OF SCIENCES
VRIJE UNIVERSITEIT
DE BOELELAAN 1081 A
1081 HV AMSTERDAM
THE NETHERLANDS
E-MAIL: aad@cs.vu.nl
E-MAIL: harry@cs.vu.nl